\documentclass[11pt]{article}
\usepackage{amssymb,amsmath,mathrsfs,amsfonts,graphicx,epsf} 
\usepackage{color}
\usepackage[dvipsnames]{xcolor}
\input xy
\xyoption{all}
\usepackage{esint}
\usepackage{fullpage}

\usepackage{authblk}

\long\def\symbolfootnote[#1]#2{\begingroup%
\def\thefootnote{\fnsymbol{footnote}}\footnote[#1]{#2}\endgroup} 

\usepackage{graphicx}
\DeclareGraphicsExtensions{.jpg}
\usepackage{color}
\usepackage{fullpage}
\usepackage[T1]{fontenc}
\newtheorem{theorem}{Theorem}
\newtheorem{corollary}[theorem]{Corollary}
\newtheorem{lemma}[theorem]{Lemma}
\newtheorem{proposition}[theorem]{Proposition}

\newtheorem{definition}[theorem]{Definition}
\newtheorem{example}[theorem]{Example}
\newtheorem{remark}[theorem]{Remark}

\newcommand{\bt}{\begin{theorem}}
\newcommand{\et}{\end{theorem}}
\newcommand{\bl}{\begin{lemma}}
\newcommand{\el}{\end{lemma}}
\newcommand{\bp}{\begin{proposition}}
\newcommand{\ep}{\end{proposition}}
\newcommand{\bc}{\begin{corollary}}
\newcommand{\ec}{\end{corollary}}
\newcommand{\bdeff}{\begin{definition}}
\newcommand{\edeff}{\end{definition}}
\newcommand{\brem}{\begin{remark}}
\newcommand{\erem}{\end{remark}}
\newcommand{\bex}{\begin{example}}
\newcommand{\eex}{\end{example}}

\newcommand{\e}{\mbox{e}}

\newcommand{\Zz}{\mathcal{Z}}
\newcommand{\Tt}{\mathcal{T}}

\newcommand{\be}{\begin{equation}}
\newcommand{\ee}{\end{equation}}

\renewcommand{\phi}{\varphi}

\newcommand{\bi}{\begin{itemize}}

\newcommand{\ei}{\end{itemize}}
\newcommand{\bd}{\begin{description}}
\newcommand{\ed}{\end{description}}

\newcommand{\bqn}{\begin{eqnarray*}}
\newcommand{\eqn}{\end{eqnarray*}}
\newcommand{\eqnn}{\nonumber\end{eqnarray}}

\newcommand{\ba}[1]{\begin{array}{#1}}
\newcommand{\ea}{\end{array}}

\newcommand{\R}{\mathbb{R}}

\newcommand{\p}{{\bar p}}

\newcommand{\VecM}{\mathrm{Vec}(M)}

\newcommand{\HH}{{\bf (H0)}}

\newcommand{\tet}{\theta}

\newcommand{\f}{\mathfrak{f}}

\newcommand{\bD}{{\mathcal D}}

\title{Cut locus and heat kernel at Grushin points of 2 dimensional almost Riemannian metrics}

\author{G.~Charlot}
\affil{Universit\'e Grenoble Alpes, 
CNRS, IF, F-38000 Grenoble, France\\
{\tt charlot@ujf-grenoble.fr}}

\begin{document}
\maketitle
\begin{abstract}
This article deals with 2d almost Riemannian structures, which are generalized Riemannian structures on manifolds of dimension 2.
Such sub-Riemannian structures can be locally defined by a pair of vector fields $(X,Y)$, playing the role of orthonormal frame, that may become
colinear on some subset. We denote $\bD=\mbox{span}(X,Y)$.
After a short introduction, I first give a description of the local cut and conjugate loci at a Grushin point $q$ 
(where $\bD_q$ has dimension 1 and $\bD_q=T_qM$) that makes appear that the cut locus may have an angle at $q$.
In a second time I describe the local cut and conjugate loci at a Riemannian point $x$ in a neighborhood of a Grushin point $q$. Finally,
applying results of \cite{BBCN}, I give the asymptotics in small time of the heat kernel $p_t(x,y)$ for $y$ 
in the same neighborhood of $q$.

\end{abstract}

\section{Introduction and definitions}

An almost Riemannian structure of dimension 2 (2-ARS for short) is a sub-Riemannian structure on a 2 dimensional manifold 
with a rank varying distribution. It is supposed to be locally defined by a pair of vector fields, playing the role of an orthonormal 
frame, that satisfies the H\"ormander condition. It defers from Riemannian geometry by the fact that the pair may become colinear.

$2$-ARSs were first studied in the context of hypoelliptic operators \cite{baouendi,FL1,grushin1}.  
They have applications to quantum control  \cite{q4,BCha,q1} and orbital transfer 
in space mechanics  \cite{BCa,tannaka}. 

The singular set $\Zz$ where the distribution $\bD$ has dimension 1 is generically a 1 dimensional embedded 
submanifold (see \cite{ABS}). There are generically three types of points : Riemannian points where the distribution 
has dimension 2; Grushin points where $\bD_q$ has dimension 1, $\bD^2_q$ has dimension 2 and $\bD_q$ 
is transverse to $\Zz$; Tangency points where $\bD_q=\bD^2_q$ has dimension 1, $\bD^3_q$ has dimension 2 and 
$\bD_q$ is tangent to $\Zz$. We denote $\Tt$ the set of tangency points.

In \cite{ABS,high-order}, a Gauss Bonnet formula is obtained for 2-ARS without tangency points. In \cite{euler}, it is 
generalized in presence of tangency points. In \cite{BCGS} a necessary and sufficient condition for two 
2-ARS to be Lipschitz equivalent is given in terms of labelled graphs associated to the structures. In \cite{camillo}
the authors show that the singular set $\Zz$ acts as a barrier for the heat flow and for quantum particules despite the fact that
geodesics can pass through $\Zz$.

A general fact in sub-Riemannian geometry is that, at points $q$ where $\bD_q\neq T_qM$, the conjugate and the cut loci of $q$
do accumulate at $q$. In \cite{BCGJ}, the local cut locus of a tangency point is described: it is an asymetric cusp, 
tangent to the distribution at $q$. In \cite{crest}, the cut locus of the wave front starting from $\Zz$ is described
in the neighborhood of a tangency point $q$: it is the union of a curve starting from $q$ and transversal to $\bD_q$ and 
of an other one tangent to $\bD_q$. 

As surprising as it may appear, nothing has been done for what concerns the local cut locus at a Grushin point $q$. 
Certainly one easily imagines that it is a $C^1$-curve transverse to $\bD_q$. In this article, I prove at the contrary that, 
for a Grushin point $q$ outside a discrete set of $\Zz-\Tt$, the cut locus of $q$ has an angle at $q$. Studying the cut and conjugate
loci for a point $q$ close to $\Zz$ and far from $\Tt$, I also give estimates of the asymptotics in small time for the heat kernel
associated with both the 2-ARS and a smooth volume.

\subsection{Definitions and basic properties}

\begin{definition}
 \label{2-ARS}
 A {\it $2$-dimensional almost-Riemannian structure} ($2$-ARS, for short) is a triple
${\mathcal S}=
(E, \f,g)$
 where:
\bi
\item $E$ is a  vector bundle of rank $2$ over a $2$ dimensional smooth manifold $M$;
\item $g$ is an Euclidean metric on $E$, that is $g_q$ is a scalar product on $E_q$ smoothly depending on $q$;
\item $\f:E\rightarrow TM$ is a morphism of vector bundles, that is $\f$ is linear from $E_q$ to $T_qM$ for any $q$.
\ei
Denoting by $\Gamma(E)$ the $C^\infty(M)$-module of smooth sections on $E$, and by 
 $\f_*:\Gamma(E)\rightarrow \VecM$ the map $\sigma\mapsto\f_*(\sigma):=\f\circ\sigma$, 
we require that the submodule of Vec$(M)$ given by $\bD=\f_*(\Gamma(E))$ is bracket generating, i.e.,
  $Lie_q(\bD)= T_qM$ for every $q\in M$.  Moreover, we require that $\f_*$ is injective.
\end{definition}

If $(\sigma_1,\sigma_2)$ is an orthonormal frame for $g$ on
 an open subset $\Omega$ of $M$, an {\it  orthonormal frame} on $\Omega$ is 
 given by $(\f_*\sigma_1,\f_*\sigma_2)$ which forms a local generator of the submodule $\bD$.  

Admissible curves, sub-Riemannian length and distance, geodesics, spheres and wave front, conjugate and cut loci, are defined 
as in the classical sub-Riemannian setting (see for example \cite{AgrBarBoscbook}).

Under the following generic\footnote{Generic means true for a residual subset of the set of morphisms $\f$ endowed with the $C^\infty$-Whitney topology} asumption \HH, only Riemannian, Grushin and tangency points can occur (see \cite{ABS}).

\bd
\item[\HH] {(i)} $\Zz$ is an embedded one-dimensional submanifold of $M$;

\hspace{.1cm}{(ii)} the points $q\in M$ where $\bD^2_q$ is
one-dimensional are isolated and at these points $\bD_q=T_q\Zz$;

\hspace{.1cm}{(iii)}  $\bD_q^3=T_qM$ for every $q\in M$.
\ed

\medskip

At Grushin points, it exists a canonical local coordinate system such that an orthonormal frame $(X,Y)$ is given by 
the normal form
$$\mbox{(NF)}\quad\quad\quad
(X=\partial_x,\quad
Y=x f(x,y)\partial_y)
$$
where $f(0,y)=1$ (see \cite{ABS,crest}). For this normal form, the natural orders
of $x$ and $y$ are respectively 1 and 2 (see \cite{bellaiche}).
The nilpotent approximation is the so called Grushin metric defined by the orthonormal frame
$$
X_{-1}=\partial_x,\quad
Y_{-1}=x\partial_y.
$$

One proves easily that there is no abnormal extremal which implies, thanks to the Pontryagin maximum principle, 
that any geodesic is the projection on $M$ of a trajectory of the Hamiltonian defined on $T^*M$ by
$$
H(\lambda,q)=\frac 1 2 ( (\lambda.X(q))^2+(\lambda.Y(q))^2)
$$
where $(X,Y)$ is an orthonormal frame of $\bD$. See 
\cite{pontryagin-book}.

In the following, we use the notation ${\mathcal A}_\ell$ for a familly of singularity. See \cite{arnold} for their definition.
A map $f$ from $\R^n$ to $\R^n$ has singularity ${\mathcal A}_2$ at $q$ if, up to a good choice of 
variables at $q$ and coordinates at $f(q)$, it can be written $(x_1,\dots,x_n)\mapsto (x_1^2,x_2,\dots,x_n)$ 
close to $x_1=\dots=x_n=0$ (a fold).
A map has singularity ${\mathcal A}_3$ at $q$ if it can be written 
$(x_1,\dots,x_n)\mapsto (x_1^3-x_1x_2,x_2,\dots,x_n)$ 
close to $x_1=\dots=x_n=0$ (a cusp).

\subsection{Results}

In this article, concerning the cut locus at a Grushin point, I prove in section \ref{section-cut}
\bt\label{cut}
Close to a Grushin point $q$ of a 2-ARS, its cut locus is the union of 
the disjoint supports of two smooth curves $\gamma_i:\;]0,\varepsilon[\;\to M$ ($i=1,2$), 
such that $\gamma_i(0)=q$, $\dot\gamma_i(0)\neq0$. 
Generically, the two vectors $\dot\gamma_1(0)$ and $\dot\gamma_2(0)$ are not colinear except at isolated points
of $\Zz\setminus\Tt$. More precisely, in the normal coordinate system such that {\rm(NF)} holds, the curves are given by
$$
\gamma_1(t)=t(-\frac 43 a,\frac \pi 2) +O(t^2), \quad\quad \gamma_2(t)=t(-\frac 43 a,-\frac \pi 2) +O(t^2).
$$
where $a=\frac{\partial f}{\partial x}(0,0)$.
\et

For what concerns the cut locus of a point close enough to the singular set $\Zz$ and far enough to the set $\Tt$, I prove in section \ref{section-cut-conjugate}
\bt\label{cut-conjugate}
Let $q_0$ be a Grushin point of a 2-ARS. If $q$ is a Riemannian point, sufficiently close to $q_0$, then the cut locus of $q$ is locally 
the union of the disjoint supports of two smooth curves $\gamma_i:\;[0,\varepsilon[\;\to M$ ($i=1,2$). The points $\gamma_i(0)$ belong
also to the first conjugate locus, they are reached by only one optimal geodesic and the corresponding singularity of the exponential map $Exp_q$ 
is of type ${\mathcal A}_3$. 
The other points of the local cut locus are reached by two optimal geodesics and do not correspond to singularities of $Exp_q$.
\et

\medskip

Considering the Riemannian volume, one can study the canonical heat equation associated with the almost Riemannian structure. 
As explained in \cite{camillo}, in the case of a compact orientable surface without tangency points, a quantum particle in such 
a structure cannot cross the singular set and the heat cannot flow through the singularity. 
This is really surprising since geodesics do cross the singular set.

Considering a smooth volume on $M$ (which is not the case of the Riemannian volume along $\Zz$), one can define
a divergence. Together with the sub-Riemannian gradient, it allows to define a Laplacian by
$$\Delta f=\mbox{div}(\nabla_g f).$$ 
Under the additional hypothesis that the manifold is complete, one gets that $\Delta$ is hypoelliptic and has a symmetric heat kernel
$p_t(x,y)$. 

In \cite{BBCN}, the authors prove in particular that if a geodesic $\gamma$ between $x$ and $y$ is such that $y$ is both in the cut locus
and the conjugate locus along $\gamma$ and if the corresponding singularity of $exp_x$ is of type ${\mathcal A}_{\ell}$ then
the contribution to the heat kernel has the following asymptotics in small time
$$
p_t(x,y)=\frac{C+O(t^{\frac{2}{\ell+1}})}{t^{\frac{n+1}{2}-\frac{1}{\ell+1}}}\e^{-d^2(x,y)/4t}.
$$
In \cite{BBN}, the authors prove that if a geodesic $\gamma$ between $x$ and $y$ is not conjugated in $y$ then
$$
p_t(x,y)=\frac{C+O(t)}{t^{\frac{n}{2}}}\e^{-d^2(x,y)/4t}.
$$
With these results, one proves easily that Theorem \ref{cut-conjugate} implies 
\bt\label{heat}
Let $q_0$ be a Grushin point of a 2-ARS. If $x$ is a Riemannian point, sufficiently close to $q_0$, 
then if $y\neq x$ is still close to $q_0$ one gets that
\bi
\item if no optimal geodesic between $x$ and $y$ is conjugated at $y$ then it exists $C$ such that
$$
p_t(x,y)=\frac{C+O(t)}{t}\e^{-d^2(x,y)/4t},
$$
\item otherwise, there is only one optimal geodesic between $x$ and $y$, which is conjugated at $y$ and it exists $C$ such that
$$
p_t(x,y)=\frac{C+O(t^{\frac12})}{t^\frac 5 4}\e^{-d^2(x,y)/4t}.
$$
\ei
\et

\section{Cut locus at a Grushin point}\label{section-cut}
Let us use the normal form (NF) at a Grushin point given by
$F_1=\left(\ba{c}1\\0\ea\right)$, $F_2=\left(\ba{c}0\\xf(x,y)\ea\right)$, with $f$ smooth such that $f(0,y)=1$.

Using the Pontryagin Maximum Principle one gets the equations
$$
\ba{rclcrcl}
\dot x & = & p_x, & &\dot p_x & = & -p_y^2 xf(x,y)(f(x,y)+x\partial_x f(x,y)),\\
&&&&&&\\
\dot y & = & p_y (xf(x,y))^2,& &\dot p_y & = &-p_y^2 x^2f(x,y)\partial_y f(x,y).
\ea
$$
Setting $\p=\frac{p_x}{p_y}$, defining the new time $s=p_y t$,
and writing $f(x,y)=1+ax+o(x,y)$, one gets the new equations
\bqn
\dot x & = & \p,\\
\dot y & = &  x^2+2ax^3+x^2o(x,y),\\
\dot \p & = & -x-3ax^2+xo(x,y).
\eqn
Initial condition is $(x=0,y=0,\p=\pm\rho)$ where $\rho=\frac {1}{p_y(0)}$.
We look at the developments of $x$, $y$ and $\p$ in the parameter $\rho$
that is
\bqn
x(\rho,s)&=&\rho x_1(s)+\rho^2 x_2(s)+O(\rho^3),\\
y(\rho,s)&=&\rho^2 y_2(s)+\rho^3 y_3(s)+O(\rho^4),\\
\p(\rho,s)&=&\rho \p_1(s)+\rho^2 \p_2(s)+O(\rho^3).
\eqn
We get the equations
$$
\ba{lll}
\dot x_1=\p_1, & \dot x_2=\p_2, & \dot y_2=x_1^2,\\
\dot \p_1=-x_1, & \dot \p_2=-3 a x_1^2 - x_2, & \dot y_3=2 a x_1^3 + 2 x_1 x_2,
\ea
$$
with the initial condition
$x_1(0)=x_2(0)=y_2(0)=y_3(0)=p_2(0)=0$ and $p_1(0)=1$.
The solution is given by
$$
\ba{lllll}
x_1(s)=\sin(s),&&\p_1(s)=\cos(s), &&y_2(s)= \frac 1 4 (2 s - \sin(2 s)),\\
&&&&\\
x_2(s)=-4 a \sin^4(\frac s 2),&&\p_2(s)=-4 a \sin^2(\frac s 2) \sin(s),& & y_3(s)=\frac {8a}{3} (1 + 2 \cos(s))\sin^4(\frac s 2).
\ea
$$
If we compute the cut locus for the nilpotent approximation (order -1), 
we find $x=0, y=\frac \pi 2 \rho^2$ for the upper part
and $x=0, y=-\frac \pi 2 \rho^2$ for the lower part.
These formulae cannot be apriori stable in the sense that the following developments could make appear terms in $\rho^2$
in the $x$ variable and hence change the "tangent at $0$" of the cut locus.

To look for the upper part of the cut locus for the normal form at order 0, we look for the cut point reached at time $t=\pi \rho_0$.
The corresponding geodesic starting with $p_1=+1$ has $\rho$ and $s$ close to $\rho_0$ and $\pi$ that is
\bqn
\rho_{^+}^{_+}&=&\rho_0+\alpha_{^+}^{_+}\rho_0^2+o(\rho_0^2),\\
s_{^+}^{_+}&=&\pi+\beta_{^+}^{_+}\rho_0+o(\rho_0).
\eqn
Since $\rho_{^+}^{_+} s_{^+}^{_+}=\rho_0 \pi$ we get immediately that $\beta_{^+}^{_+}=-\alpha_{^+}^{_+}\pi$ and
\bqn
x_{^+}^{_+}&=&(\alpha_{^+}^{_+}\pi-4a)\rho_0^2+O(\rho_0^3),\\
y_{^+}^{_+}&=&\frac \pi 2 \rho_0^2+(\alpha_{^+}^{_+}\pi-\frac{8a}{3})\rho_0^3+O(\rho_0^4).
\eqn
The corresponding geodesic starting with $p_1=-1$ has $\rho$ and $s$ 
\bqn
\rho_{^+}^{_-}&=&\rho_0+\alpha_{^+}^{_-}\rho_0^2+o(\rho_0^2),\\
s_{^+}^{_-}&=&\pi+\beta_{^+}^{_-}\rho_0+o(\rho_0),
\eqn
Since $\rho_{^+}^{_-} s_{^+}^{_-}=\rho_0 \pi$ we get immediately that $\beta_{^+}^{_-}=-\alpha_{^+}^{_-}\pi$ and
\bqn
x_{^+}^{_-}&=&(-\alpha_{^+}^{_-}\pi-4a)\rho_0^2+O(\rho_0^3),\\
y_{^+}^{_-}&=&\frac \pi 2 \rho_0^2+(\alpha_{^+}^{_-}\pi+\frac{8a}{3})\rho_0^3+O(\rho_0^4).
\eqn
Since $(x_{^+}^{_+},y_{^+}^{_+})$ and $(x_{^+}^{_-},y_{^+}^{_-})$ should be the same point we get
\bqn
\alpha_{^+}^{_+}\pi-4a&=&-\alpha_{^+}^{_-}\pi-4a,\\
\alpha_{^+}^{_+}\pi-\frac{8a}{3}&=&\alpha_{^+}^{_-}\pi+\frac{8a}{3},
\eqn
which implies $\alpha_{^+}^{_+}=-\alpha_{^+}^{_-}=\frac{8a}{3\pi}$ and that the cut point is
$$
(x_{cut}^+, y_{cut}^+)=\rho_0^2(-\frac 43 a,\frac \pi 2) +O(\rho_0^3).
$$

The same computation for the lower part of the cut locus gives
$$
(x_{cut}^-,y_{cut}^-)=\rho_0^2(-\frac 43 a,-\frac \pi 2)+O(\rho_0^3).
$$
The two formulae for the upper and lower parts are stable in the sense that the terms that could be added by further
developments would be of order at least 3 in $\rho_0$ and hence would not change the tangent at 0.
We can now conclude that the cut locus has a corner at 0 when $a\neq 0$ and none when $a=0$. Theorem \ref{cut} is proved.

\section{Singularities of the exponential map from a point $q$ close enough
to a given Grushin point $q_0$. Applications to the heat kernel at $q$}\label{section-cut-conjugate}

\subsection{Study in the Grushin plane $(\partial_x,x\partial_y)$.}

In this section I discribe the cut and conjugate loci at the point $(-1,0)$ in the Grushin plane.
In the case of the Grushin plane, whose orthonormal frame is given by $(\partial_x,x\partial_y)$, one can compute explicitely the geodesics.
It was done in \cite{camillo}, and the geodesics from $(-1,0)$ are given by
\bqn
x(\tet,t)&=&-\frac{\sin(\tet - t\sin(\tet))}{\sin(\tet)},\\
y(\tet,t)&=&\frac{2 t\sin(\tet) - 2 \cos(\tet)\sin(\tet) + \sin(2 \tet - 2 t \sin(\tet))}{4 \sin^2(\tet)}.
\eqn
if $\tet\neq 0[\pi]$ or
\bqn
x(\tet,t)&=&-1+(-1)^{\tet/\pi} t,\\
y(\tet,t)&=&0,
\eqn
if $\tet=0[\pi]$. We denote $\gamma(\tet,t)=(x(\tet,t),y(\tet,t))$. The geodesics parameterized by arclength are the $t\mapsto\gamma(\tet,t)$.

\medskip

\noindent {\bf The cut locus.}

For $\tet\neq0[\pi]$ : 
$$
\gamma\left(\theta, \frac \pi{|sin(\tet)|}\right)=\gamma\left(\pi-\theta, \frac \pi{|\sin(\tet)}|\right)=\left(1,\frac{\pi}{2\sin(\tet)|\sin(\tet)|}\right).
$$
These geodesics being normal geodesics, it implies that
their cut time is less or equal to $\frac \pi{|\sin(\tet)|}$. Moreover, it is very easy to prove that any geodesic corresponding to $\tet=0[\pi]$ is optimal until any time $t$. Let us prove that the cut time is in fact $\frac \pi{\sin(\tet)}$ if $\tet\neq0[\pi]$.

The sphere at time $t>\pi$ is contained in 
$\{\gamma(\tet,t)\; | \; \tet\in \Theta_t\}$ where $$\Theta_t=[-\arcsin(\frac\pi t),\arcsin(\frac\pi t)]\cup[\pi-\arcsin(\frac\pi t),\pi+\arcsin(\frac\pi t)]$$ 
since a geodesic with $t\geq\frac\pi{|\sin(\tet)|}$ is no more optimal. We are going to prove that $\tet\mapsto x(\tet,t)$ is strictly decreasing on 
$[0,\arcsin(\frac\pi t)]$ and on $[\pi-\arcsin(\frac\pi t),\pi]$, which correspond to initial conditions constructing the upper
part of the synthesis. If one add the fact the $\gamma(\arcsin(\frac\pi t),t)=\gamma(\pi-\arcsin(\frac\pi t),t)$ it will have proved that the upper part of the sphere, being connex, is the set $\{\gamma(\tet,t)\; | \; \tet\in \Theta_t\}$. Which implies that the cut locus of $(-1,0)$ is attained at time $\frac \pi{\sin(\tet)}$ and hence that the cut locus is
$$
\{(x,y)\; | \; x=1 \mbox{ and } |y|\geq \pi/2\}.
$$
In order to finish the proof, let us compute for $\frac{\partial x}{\partial \tet}(\tet,t)$ and prove that for $\tet\in]0,\pi[$ and 
$t\in[0,\frac\pi{\sin(\tet)}]$ it is not positive.
\begin{eqnarray*}
-\sin^2(\tet)\frac{\partial x}{\partial \tet}(\tet,t)&=&\sin(t\sin(\tet))-t\cos(\tet)\sin(\tet)\cos(\tet-t\sin(\tet))\\
&=&\sin(u)-u\cos(\tet)\cos(\tet-u)\\
&=&\sin(u)(1-u\cos(\tet)\sin(\tet))-\cos(u) u \cos^2(\tet)
\end{eqnarray*}
where $u=t\sin(\tet)$. For $t=0$, $\frac{\partial x}{\partial \tet}(\tet,t)=0$.
And for $t$ small $\frac{\partial x}{\partial \tet}(\tet,t)\sim -t\sin(\tet)<0$. Now we can conclude by proving that 
$\frac{\partial x}{\partial \tet}(\tet,t)\neq 0$ for $0<t<\frac\pi{\sin(\tet)}$.

With the last equation we can see that $\frac{\partial x}{\partial \tet}(\tet,t)= 0$ if and only if $(\cos(u),\sin(u))$ is parallel to 
$(1-u\cos(\tet)\sin(\tet),u\cos^2(\tet))$. We just have to prove now that with $\tet\in]0,\pi[$ and $u\in]0,\pi[$ it is not possible.
This is a relatively simple exercice of geometry. Let us make the proof for $\tet\in]0,\pi/2[$, the proof being the same for $\tet\in]\pi/2,\pi[$
and very easy for $\tet=\pi/2$. 

Let us fix the angles of the two vectors to be zero for $u=0$.
The first vector, $(\cos(u),\sin(u))$ has for angle $u$. The second vector $(1-u\cos(\tet)\sin(\tet),u\cos^2(\tet))=(1,0)-u\cos(\tet)(\sin(\tet),-\cos(\tet))$
has norm larger than $\cos(\tet)$ if $u\neq \sin(\tet)$ and a derivative with respect to $u$ of norm $\cos(\tet)$. Hence the derivative of its angle 
with respect to $u$ is less than 1 (and positive) for $u\neq \sin(\tet)$ which allows to prove that its angle is positive and less than $u$ for $u>0$. 
As a conclusion, if $u<\pi$, the angles of the two vectors cannot be equal modulo $\pi$. Which finishes the proof.
\brem
One computes easily that $\frac{\partial\gamma}{\partial t}(\tet,\frac{\pi}{\sin(\tet)})$
and $\frac{\partial\gamma}{\partial t}(\pi-\tet,\frac{\pi}{\sin(\tet)})$ are not parallel hence the wave front
is transversal to itself along the cut locus.
\erem

\medskip
\noindent{\bf The first conjugate locus.}
The Jacobian of the map $\gamma$ is
$$
Jac(\tet,t)=\frac{t \cos (\tet) \cos (\tet-t \sin (\tet)) \sin (\tet)-\sin (t \sin (\tet))}{\sin^3(\tet)}.
$$
when $\tet\neq0[\pi]$. One proves easily that if $\tet\neq 0[\pi]$, then there is a conjugate time $t_\tet$.
Moreover at $t=t_\tet$ one has that the vectors $(\cos(t_\tet\sin(\tet)),\sin(t_\tet\sin(\tet)))$ and $(1-t_\tet\cos(\tet)\sin^2(\tet),t_\tet\cos^2(\tet)\sin(\tet))$ are parallel since 
$$
0=Jac(\tet,t_\tet)=\frac{t_\tet\cos^2(\tet)\sin(\tet)\cos(t_\tet\sin(\tet)-(1-t_\tet\cos(\tet)\sin^2(\tet))\sin(t_\tet\sin(\tet)))}{\sin^3(\tet)}.
$$
One can compute
$$
\frac{\partial\gamma}{\partial\tet}(\tet,t)=\frac{1}{\sin^3(\tet)}Jac(\tet,t) \left(\sin(\tet),-\cos(\tet - t \sin(\tet))\right),
$$
which proves that at the conjugate time $\frac{\partial\gamma}{\partial\tet}(\tet,t_\tet)=0$.

In order to understand which singularity has the map $\gamma$ at the conjugate time, let compute
\begin{eqnarray*}
\frac{\partial^2\gamma}{\partial\tet^2}(\tet,t_\tet)&=&\frac{1}{\sin^3(\tet)}
\frac{\partial Jac}{\partial\tet}(\tet,t_\tet)\left(\sin(\tet),-\cos(\tet -  t_\tet\sin(\tet))\right),\\
\frac{\partial\gamma}{\partial t}(\tet,t_\tet)&=&\frac{1}{\sin(\tet)}(\cos(\tet - t_\tet \sin(\tet))\sin(\tet), \sin(\tet - t_\tet \sin(\tet))^2).
\end{eqnarray*}
Now one can compute the determinant of these two vectors and find at $t=t_\tet$
\begin{eqnarray*}
Det(\tet)&=&\frac{1}{\sin^4(\tet)}\frac{\partial Jac}{\partial\tet}(\tet,t_\tet)
\left|
\ba{cc}
\sin(\tet) & -\cos(\tet -  t_\tet\sin(\tet))\\
\cos(\tet - t \sin(\tet))\sin(\tet)& \sin(\tet - t \sin(\tet))^2
\ea
\right|,\\
&=&\frac{1}{\sin^3(\tet)}\frac{\partial Jac}{\partial\tet}(\tet,t_\tet).
\end{eqnarray*}
As a consequence, one can conclude that the singularity is of type ${\mathcal A}_2$ if and only if $\frac{\partial Jac}{\partial\tet}(\tet,t_\tet)\neq0$. But one can compute that
\begin{eqnarray*}
\frac{\partial Jac}{\partial\tet}(\tet,t_\tet)&=&
\frac{t_\tet\sin(\tet)}{4}
\left(
\sin(t_\tet\sin(\tet))
(2-3t_\tet \cos(\tet)+6\cos(2\tet)-t_\tet\cos(3\tet))\right.\\
&&\left.+\cos(t_\tet\sin(\tet))
(t_\tet\sin(\tet)-6\sin(2\tet)+t_\tet\sin(2\tet))
\right).
\end{eqnarray*}
which implies that $\frac{\partial Jac}{\partial\tet}(\tet,t_\tet)=0$ if and only if
$$
(t_\tet\cos^2(\tet)\sin(\tet))
(2-3t_\tet \cos(\tet)+6\cos(2\tet)-t_\tet\cos(3\tet))\hspace{3cm}\
$$
$$\ \hspace{3cm}
+(1-t_\tet\cos(\tet)\sin^2(\tet))
(t_\tet\sin(\tet)-6\sin(2\tet)+t_\tet\sin(2\tet))=0
$$
Simplifying this last expression one finds 
$$
-(6 + t_\tet^2 - 6 t_\tet \cos(\tet) + t_\tet^2 \cos(2\tet)) \sin(2 \tet)=0
$$
One proves easily that this expression is zero if and only if $\tet=0[\frac \pi 2]$. Hence, for any value of $\tet\neq 0[\frac \pi 2]$, the singularity at the conjugate time is of type ${\mathcal A}_2$. 

To understand completely the optimal synthesis from $(-1,0)$, it remains to understand which singularity is for $\tet=\frac \pi 2[\pi]$ at its conjugate time, which is equal to its cut-time $t_\frac{\pi}{2}=\pi$. In order to do that
let us make the change of variables $\tet=\frac\pi 2+\tet_1$, $s=\pi+s_1$ and still denote the map $\gamma$. Then the Taylor series up to order 3 of $\gamma$ at $(\tet_1=0,s_1=0)$ is
$$
\gamma(\tet_1,s_1)=(1+-\frac{s_1^2}{2}+\tet_1 s_1+
\frac{1}{2} \left(\pi  s_1 \tet_1^2-\pi \tet_1^3\right),
\frac{\pi }{2}+s_1-\frac{s_1^3}{3}+\tet_1 s_1^2-\frac{\tet_1^2 s_1}{2}
)+o_3(\tet_1,s_1).
$$
Making the change of coordinates $x_1=x-1+\frac{(y-\frac\pi 2)^2}{2}$, $y_1=y-\frac\pi 2$, still
denoting the map $\gamma$, we found
$$
\gamma(\tet_1,s_1)=(\tet_1 s_1+
\frac{1}{2} \left(\pi  s_1 \tet_1^2-\pi \tet_1^3\right),
s_1-\frac{s_1^3}{3}+\tet_1 s_1^2-\frac{\tet_1^2 s_1}{2})+o_3(\tet_1,s_1).
$$
By changing for the variables $\tet_2=\tet_1 +
\frac{1}{2} \pi  \tet_1^2$, $s_2=s_1-\frac{s_1^3}{3}+\tet_1 s_1^2-\frac{\tet_1^2 s_1}{2}$, still
denoting the map $\gamma$, we found
$$
\gamma(\tet_2,s_2)=(\tet_2 s_2-\frac{1}{2} \pi \tet_2^3,s_2)+o_3(\tet_2,s_2).
$$
Finally, making the changes $y_2=-\frac \pi 2 y_1$, $x_2=-\frac \pi 2 x_1$,
$s_3=\frac \pi 2 s_2$, $\tet_3=\tet_2$, we get
$$
\gamma(\tet_3,s_3)=(\tet_3^3-\tet_3 s_3,s_3)+o_3(\tet_3,s_3),
$$
which proves that the singularity at the first conjugate locus for $\tet=\frac\pi 2[\pi]$ is ${\mathcal A}_3$.

\subsection{Local cut and conjugate loci of a Riemannian point lying in the neighborhood of a Grushin point}

Assume that we are close enough to a Grushin point in order
the normal form $(F_1,F_2)$ applies. In that case we can compute the
jets of the exponential map from $(-a,0)$ with respect to the
small parameter $a$. From the computation done in the case of the nilpotent approximation,
we can deduce that $\gamma_a$ the exponential map from $(-a,0)$ has the following expression
\bqn
x_a(\tet,t)&=&-a\left(\frac{\sin(\tet - \frac{t}{a}\sin(\tet))}{\sin(\tet)}\right)+O(a^2),\\
y_a(\tet,t)&=&a^2\left(\frac{2\frac{t}{a}\sin(\tet) - 2 \cos(\tet)\sin(\tet) + \sin(2 \tet - 2 \frac{t}{a} \sin(\tet))}{4 \sin^2(\tet)}\right)+O(a^3),
\eqn
which proves that $(\theta,s)\mapsto(\frac{x_a}{a},\frac{y_a}{a^2})(\tet,as)=(x,y)(\tet,s)+O(a)$.
As seen before, the map $\gamma$ has only stable singularities at its first conjugate locus
and the front is transversal to itself at the cut locus outside the conjugate points. It implies
that $(\theta,s)\mapsto(\frac{x_a}{a},\frac{y_a}{a^2})(\tet,as)$, as one parameter familly of perturbation
of the map $(x,y)$, has the same singularities at the conjugate locus and have a wave front transversal
to itself at the cut locus outside the conjugate locus, at least on any compact for $a$ small enough.
We can assume that the compact corresponds to $0\leq s\leq 2\pi$, that way we are sure that the 
cut-conjugate points are in the interior of the compact set.
Hence the map $(\theta,t)\mapsto(\frac{x_a}{a},\frac{y_a}{a^2})(\tet,t)$  has the same singularities at 
the conjugate locus and have a wave front transversal to itself at the cut locus outside the conjugate locus,
for $a$ small enough and $0\leq t\leq 2\pi a$. Which implies that $\gamma_a$ has for local cut locus
and local conjugate locus the same picture as the map $\gamma$: only two cut-conjugate points, where
the singularity is ${\mathcal A}_3$, other first conjugate points with singularity ${\mathcal A}_2$ and
the cut locus is the union of
two disjoint curves issued  from the cut-conjugate points along which the wave front is transversal to itself.
This finishes the proof of theorem \ref{cut-conjugate}.

\bibliography{biblio-grushin2D}

\bibliographystyle{plain}

\end{document}